\documentclass[11pt, reqno]{amsart}

\usepackage[all]{xy}
\usepackage{amssymb}

\newtheorem{thm}{Theorem}[section]
\newtheorem{prop}[thm]{Proposition}
\newtheorem{lem}[thm]{Lemma}
\newtheorem{cor}[thm]{Corollary}

\theoremstyle{definition}

\newtheorem{remark}[thm]{Remark}

\newtheorem{reduction}[thm]{Reduction}
\newtheorem*{ack}{Acknowledgment}

\numberwithin{equation}{section}

\newcommand{\UD}{\operatorname{UD}}

\newcommand{\Spec}{\operatorname{Spec}}

\newcommand{\Hom}{\operatorname{Hom}}

\newcommand{\bbG}{\mathbb{G}}

\newcommand{\bbZ}{\mathbb{Z}}

\newcommand{\bbQ}{\mathbb{Q}}

\newcommand{\cV}{\mathcal{V}}

\newcommand{\cG}{\mathcal{G}}
\newcommand{\cH}{\mathcal{H}}

\newcommand{\isomoto}{\overset{\sim}{\to}}

\newcommand{\brokrarr}{\vphantom{\to}\mathrel{\smash{{-}{\rightarrow}}}}

\newcommand{\lf}{\mathopen}

\let\r=\mathclose

\newcommand{\Ker}{\operatorname{Ker}}
\newcommand{\Aut}{\operatorname{Aut}}
\newcommand{\sdp}{\mathbin{{>}\!{\triangleleft}}}

\newcommand{\GL}{{\operatorname{GL}}}

\newcommand{\PGL}{{\operatorname{PGL}}}
\newcommand{\PGLn}{{\operatorname{PGL}_n}}
\newcommand{\rank}{\operatorname{rank}}
\newcommand{\Stab}{\operatorname{Stab}}

\newcommand{\diag}{\operatorname{diag}}
\newcommand{\ed}{\operatorname{ed}}
\newcommand{\Gal}{\operatorname{Gal}}
\newcommand{\Galois}{\Gal}
\newcommand{\lra}{\longrightarrow}
\newcommand{\Mat}{{\operatorname{M}}}
\newcommand{\Mn}{\Mat_n}

\newcommand{\trdeg}{\operatorname{trdeg}}
\newcommand{\Sym}{{\operatorname{\mathcal{S}}}}

\newcommand{\wed}[1]{\bigwedge^2{#1}}
\newcommand{\sym}[1]{\operatorname{\mathsf{Sym}}^2{#1}}
\newcommand{\ind}[2]{\!\!\uparrow_{#1}^{#2}}
\def\onto{\twoheadrightarrow}
\newcommand{\Ext}{\operatorname{Ext}}
\newcommand{\Res}{\operatorname{Res}}

\newdir{ _{(}}{{}*!/-3pt/\dir_{(}}

\let\to=\longrightarrow

\begin{document}
\title[Lattices and division algebras]{Lattices and 
parameter reduction in division algebras}

\author{M. Lorenz}
\address{Department of Mathematics, Temple University,
	Philadelphia, PA 19122-6094}
\email{lorenz@math.temple.edu}
\thanks{M. Lorenz was supported in part by NSF grant DMS-9618521}

\author{Z. Reichstein}
\address{Department of Mathematics, Oregon State University,
	Corvallis, OR 97331}
\email{zinovy@math.orst.edu}
\thanks{Z. Reichstein was supported in part by NSF grant DMS-9801675}


\subjclass{16K20, 14L30, 20C10, 20J06}
\keywords{division algebra, integral representation, lattice, algebraic transformation group,
essential dimension, crossed product, central simple algebra}

\begin{abstract}
Let $k$ be an algebraically closed field of characteristic $0$ and
let $D$ be a division algebra whose center $F$ contains $k$.
We shall say that $D$ can be reduced to $r$ parameters if
we can write $D \simeq D_0 \otimes_{F_0} F$, where $D_0$ is
a division algebra, the center $F_0$ of $D_0$ 
contains $k$ and $\trdeg_k(F_0) = r$.  

We show that every division algebra of odd degree $n \geq 5$ can be reduced to
$\leq {\frac{1}{2}}(n-1)(n-2)$ parameters. Moreover, every
crossed product division algebra of degree $n \geq 4$ can be reduced to
$ \leq (\lfloor \log_2(n) \rfloor - 1)n + 1$ parameters.
Our proofs of these results rely on lattice-theoretic techniques.
\end{abstract}

%
%

\maketitle
\tableofcontents


\section{Introduction}

Throughout this paper $k$ denotes a (fixed) algebraically closed base
field of characteristic zero.
Let $K$ be a field containing $k$ and let
$A$ be a finite-dimensional $K$-algebra.
We would like to write  $A$ as
$A = A_0 \otimes_{K_0}K$ for some $K_0$-algebra $A_0$ over an intermediate field
$k\subseteq K_0\subseteq K$ with $\trdeg_k(K_0)$ as low as possible;
the minimal value of $\trdeg_k(K_0)$ will be denoted by $\tau(A)$.
Note that if $\trdeg_k(K_0) < \trdeg_k(K)$ then passing from
$A$ to $A_0$ may be viewed as ``parameter reduction" in $A$.

We shall be particularly interested in the case where $A = \UD(n)$ is
the universal division algebra of degree $n$ and $K$ is the center of $\UD(n)$
which we shall denote by $Z(n)$. Recall that $\UD(n)$ is
the subalgebra of $\Mn(k(x_{ij}, y_{ij}))$ generated
(as a division algebra) by two generic $n \times n$-matrices
$X = (x_{ij})$ and $Y = (y_{ij})$,
where $x_{ij}$ and $y_{ij}$ are $2n^2$ independent variables over $k$;
see, e.g.,~\cite[Section II.1]{procesi} or~\cite[Section 3.2]{rowen-pi}.
We will denote $\tau(\UD(n))$ by
$d(n)$.  It is easy to show that $d(n) \geq \tau(A)$ for any central
simple algebra $A$ of degree $n$ whose center contains $k$
(see, e.g., \cite[Lemma 9.2]{reichstein1});
in other words, every
central simple algebra of degree $n$ can be ``reduced to at most $d(n)$
parameters". In the language of~\cite{reichstein1}, $d(n) = \ed(\PGLn)$,
where $\ed$ denotes the essential dimension; see \cite[Lemma 9.2]{reichstein1}.

To the best of our knowledge, the earliest attempt
to determine the value of $d(n)$ is due to
Procesi, who showed that $d(n) \leq n^2$; see~\cite[Thm. 2.1]{procesi}.
Note that if $\UD(n)$ is cyclic then $d(n) = 2$,
because we can take $A_0$ to be a symbol algebra;
cf.~\cite[Lemma 9.4]{reichstein1}. This is
known to be the case for $n = 2$, $3$ and $6$. For other $n$
the exact value of $d(n)$ is not known. However, the
following inequalities hold:
\begin{gather} 
  d(n)  \leq  n^2 -2n \quad 
 \text{\rm (\cite[Proposition 4.5]{reichstein1})} \, ,  \label{e1.1} \\
 d(n)  \leq  d(nm) \leq d(n) + d(m) \quad
 \text{\rm if $(n, m) = 1$} \quad 
 \text{\rm (\cite[Section 9.4]{reichstein1})} \, , \label{e1.2} \\ 
 d(n^r)  \geq   2r \quad 
 \text{\rm (\cite[Theorem 16.1]{reichstein2})} \, , \label{e1.3} \\
 d(n)  \leq  \tfrac{1}{2}(n-1)(n-2) + n  \quad \text{\rm if $n$ is odd} 
 \quad
 \text{\rm (\cite{rowen-brauer}; cf.~\cite[Section 9.3]{reichstein1})} \, . \label{e1.4}
\end{gather}
The last inequality is due to Rowen. In this paper we will sharpen it
by showing that, in fact, $d(n) \leq \frac{1}{2}(n-1)(n-2)$ for
every odd $n \geq 5$. Moreover, in $\UD(n)$, reduction to this number
of parameters can be arranged in a particularly nice fashion:

\begin{thm} \label{thm1} Let $n \geq 5$ be an odd integer, $\UD(n)$ be
the universal division algebra of degree $n$ and
$Z(n)$ be its center.  Then there exists a subfield $F$ of $Z(n)$ and
a division algebra $D$ of degree $n$ with center $F$ such that
\begin{enumerate}
\item
$\UD(n) = D \otimes_F Z(n)$,
\item
$\trdeg_k(F) = \frac{1}{2}(n-1)(n-2)$ and
\item
$Z(n)$ is a rational extension of $F$.
\end{enumerate}
In particular, $d(n) = \ed(\PGLn) \leq \frac{1}{2}(n-1)(n-2)$.
\end{thm}

In the course of the proof of Theorem~\ref{thm1} we will obtain
an explicit description of the center $F$ of $D$:
$F \simeq k(\wed{A_{n-1}})^{\Sym_n}$, where
$\Sym_n$ denotes the symmetric group on $n$ symbols and
\begin{equation} \label{e.A}
A_{n-1} = \{ (a_1, \dots, a_n)
\in \bbZ^n  \, | \, a_1 + \dots + a_n = 0 \} \, ,
\end{equation}
with the natural $\Sym_n$-action. Our argument relies on the results of
\cite{ll}, where the symmetric square $\sym{A_{n-1}}$ is shown to be
stably permutation for $n$ odd; see Proposition~\ref{P:ll} below.

For our next result, recall that if $A$ is a central simple algebra
of degree $n$ with center $F$ and $L$ is a subfield of $A$ then
$L$ is called {\em strictly maximal} if $F \subset L$ and $[L:F] = n$.

\begin{thm} \label{thm2} Let $A$ be a finite-dimensional
central simple algebra of degree $n$ with center $F$, 
$L$ be a strictly maximal subfield
of $A$, $L^{norm}$ be the normal closure of $L$ over $F$, 
and $\cG = \Galois(L^{norm}/F)$.
Suppose $\cG$ is generated by $r$ elements 
together with $\cH=\Galois(L^{norm}/L)$.
If either $r\ge 2$ or $\cH \neq \{ 1\}$ then $\tau(A) \leq r|\cG| - n + 1$.
\end{thm}

If $A$ be a central simple algebra of degree $n$ then
the upper bounds we have for $\tau(A)$ (or, equivalently, for $d(n)$),
are all quadratic in $n$; see~\eqref{e1.1},~\eqref{e1.4} 
and Theorem~\ref{thm1}. 
However, if we assume that $A$ is a crossed product, 
Theorem~\ref{thm2} yields an asymptotically better bound: 

\begin{cor} \label{cor2a}
Suppose a group $\cG$ of order $n$ can be generated by $r\ge 2$ elements. 
Then $\tau(A) \leq (r-1)n + 1$ for any $\cG$-crossed product central
simple algebra $A$. In particular, 
$\tau(A) \leq (\lfloor \log_2(n) \rfloor - 1)n + 1$, for 
any crossed product central simple algebra of degree $n \geq 4$.
Here, as usual, $\lfloor x \rfloor$ denotes the largest integer $\leq x$.
\qed
\end{cor}

Recall that $A$ is called a $\cG$-{\em crossed product} if it contains 
a strictly maximal subfield $L$, such that $L/F$ is a Galois extension 
and $\Galois(L/F) = \cG$; cf.~\cite[Definition 3.1.23]{rowen-pi}. Thus 
the first assertion of the corollary is an immediate consequence 
of Theorem~\ref{thm2}.  The second assertion follows from the first,
because any group of order $n$ can be generated 
by $r  \leq \log_2(n)$ elements. 
(Indeed, $|\lf< \cG_0, g \r>| \geq 2 |\cG_0|$ 
for any subgroup $\cG_0$ of $\cG$ and any $g \in \cG \setminus \cG_0$.)
Note also that 
$\lfloor \log_2(n) \rfloor \geq 2$ for any $n \ge 4$. 

The case of central simple algebras of degree 4 is of special interest
since, by a theorem of Albert, every such algebra
is a $\bbZ/2\bbZ \times \bbZ/2\bbZ$-crossed product;
see e.g.,~\cite[Theorem 3.2.28]{rowen-pi}.
Thus Corollary~\ref{cor2a} says that $d(4) = \tau(\UD(4)) \leq 5$.
On the other hand, $d(4) \geq 4$ by~\eqref{e1.3}. This proves:
  
\begin{cor} \label{cor3}
$d(4) = \ed(\PGL_4) = 4$ or $5$.
\qed
\end{cor}

The rest of this paper is structured as follows. In Section~\ref{sect2}
we discuss preliminary material from invariant theory 
and the theory of $\cG$-lattices.  In Section~\ref{sect2.5} we explain 
how $\cG$-lattices can be used to give an upper bound
on essential dimensions of certain groups. We prove 
Theorem~\ref{thm1} in Section~\ref{sect3} and Theorem~\ref{thm2}  
in Section~\ref{sect5}. In Section~\ref{sect6} we show that
the methods of this paper cannot be used to decide 
whether the exact value of $d(4)$ equals $4$ or $5$.

\section{Preliminaries}
\label{sect2}

\subsection{$G$-varieties}
A $G$-{\em variety} $X$ is an algebraic variety with a (regular) action of
an algebraic group $G$. If $G$ acts freely (i.e., with trivial stabilizers)
on a dense open subset of $X$, then $X$ is called a {\em generically free}
$G$-variety. 

A dominant rational map $\pi \colon X \brokrarr Y$ is called the {\em rational
quotient map} if $k(Y) = k(X)^G$ and $\pi^{*} \colon k(Y) \hookrightarrow k(X)$
is the natural inclusion $k(X)^G \hookrightarrow k(X)$. 
We will usually denote the rational quotient $Y$ by $X/G$; note
that $X/G$ is only defined up to birational equivalence. 
By a theorem of Rosenlicht a rational quotient map separates points 
in general position in $X$; see~\cite[Theorem 2]{rosenlicht1} 
and~\cite{rosenlicht2} (also cf.~\cite[Theorem 2.3]{pv}).
In other words, there exists a dense open subset $U$ of $X$ such that
$f$ is regular on $U$ and $x, y \in U$ lie in the same $G$-orbit iff
$f(x) = f(y)$. 

We will not generally assume that $X$ is irreducible; however, we will
always require $X$ to be {\em primitive}. 
This means that $G$ transitively permutes
the irreducible components of $X$; equivalently,
$X/G$ is irreducible, i.e., $k(X)^G$ is a field; 
cf.~\cite[Section 2.2]{reichstein1}.  Note that an irreducible $G$-variety
is always primitive, and, if $G$ is connected, a primitive variety
is necessarily irreducible. Thus the notion of a primitive variety is only 
of interest if the group $G$ is disconnected. 

If $N$ is a normal subgroup of $G$ then the $G$-action on $X$
induces a (rational) $G/N$-action on $X/N$; moreover, one can choose
a model $Y$ of $X/N$ such that the $G/N$-action on $Y$
is regular; see~\cite[Proposition 2.6]{pv},~\cite[Remark 2.6]{reichstein1}.

\begin{lem} \label{lem2.1}
Let $N$ be a normal subgroup of $G$ and let $X$ be a $G$-variety. Then
$X$ is generically free as a $G$-variety if and only if
\begin{enumerate}
\item
$X$ is generically free as an $N$-variety and
\item
$X/N$ is generically free as a $G/N$-variety.
\end{enumerate}
\end{lem}

\begin{proof} 
Assume (a) and (b) hold. Choose $x$ in $X$ in general position
and suppose $g \in \Stab(x)$.  Then (b) implies that $g\in N$ 
and (a) says that $g = 1$.  This shows that the $G$-action on $X$
is generically free.  The converse in proved in a similar manner.
\end{proof}

\subsection{$(G,H)$-sections and compressions}
\label{sect2.2}

Let $X$ be a $G$-variety and let $\pi \colon X \brokrarr X/G$ be the
rational quotient map. Furthermore, let $H$ be a closed subgroup of $G$.
An $H$-invariant subvariety $S$ of $X$
is called a {\em $(G,H)$-section} if the following conditions
are satisfied.

\medskip
(i) $\pi(S)$ contains an open
dense subset of $X/G$ and

\smallskip
(ii) if $x$ is a point in general position
in $S$ then $gx \in S$ if and only if $g \in H$.

\medskip\noindent
Recall that a {\em $G$-compression} $X \brokrarr Y$ is a dominant rational
map of generically free $G$-varieties.

\begin{lem} \label{lem2.2}
Suppose $S$ is a $(G,H)$-section of $X$. Then
\begin{enumerate}
\item 
$k(S)^H = k(X)^G$ and
\item 
any $H$-compression $S \brokrarr S'$ lifts to a $G$-compression
$X \brokrarr X'$, where $S'$ is a $(G,H)$-section of $X'$.
\end{enumerate}
\end{lem}

\begin{proof}
(a) See~\cite[Section 2.8]{pv},~\cite[1.7.2]{popov}
or~\cite[Lemma 2.11]{reichstein1}.

(b) Note that $X$ and $G \ast_H S$  are birationally equivalent as $G$-varieties,
where $G \ast_H S$ is defined as $G \times S / H$ for the $H$-action
given by $h(g, s) = (gh^{-1}, hs)$; see~\cite[Theorem 1.7.5]{popov}
or~\cite[Lemma 2.14]{reichstein1}.
Now we set $X' = G \ast_H S'$
and extend $f$ to a rational map $X \brokrarr X'$ by
$f(g, s) \lra (g, f(s))$. It is easy to see that
this map has the desired properties; cf. the proof
of~\cite[Lemma 4.1]{reichstein1}.
\end{proof}

\subsection{$\cG$-lattices}
Let $\cG$ be a finite group.
A \emph{$\cG$-lattice} is a (left) module over the integral
group ring $\bbZ[\cG]$ that is free of finite rank as a $\bbZ$-module.
A $\cG$-lattice $M$ is called
\begin{itemize}
\item \emph{faithful} if the structure map $\cG\to\Aut_{\bbZ}(M)$ is injective, and
\item a \emph{permutation lattice} if $M$ has a $\bbZ$-basis that
is permuted by $\cG$.
\end{itemize}
For any $\cG$-lattice $M$, the $\cG$-action on $M$ extends canonically to actions of 
$\cG$ on the group algebra $k[M]$ and on the field of fractions $k(M)$ of $k[M]$.
The operative fact concerning the $\cG$-fields (i.e., fields with
$\cG$-action) of the form
$k(M)$ for our purposes is the following result of Masuda~\cite{masuda}; cf.~also
\cite[proof of Prop. (1.5)]{lenstra}.

\begin{prop}\label{P:lenstra}
Let $0\to M \to E \to P \to 0$ be an exact sequence of $\cG$-lattices with
$M$ faithful and $P$ permutation. Then, as $\cG$-fields,
$k(E)\simeq k(M)(t_1,\ldots,t_r)$, where the elements $t_i$ are
$\cG$-invariant and transcendental over $k(M)$  and $r=\rank P$.
\end{prop}

\subsection{The symmetric and exterior squares}
Let $M$ be a $\cG$-lattice.
By definition, the symmetric square $\sym{M}$ is the quotient of
$M^{\otimes 2}=M\otimes M$ modulo the subgroup generated by the
elements $m\otimes m' - m'\otimes m$ for $m,m'\in M$.
Similarly, the exterior square $\wed{M}$
is the quotient of $M^{\otimes 2}$ modulo the subgroup generated by the
elements $m\otimes m$, as $m$ ranges over $M$.
The action of $\cG$ on $M^{\otimes 2}$ restricts down 
to $\sym{M}$ and $\wed{M}$, making each a $\cG$-lattice.
The $\cG$-lattice $\wed{M}$ can be
identified with the sublattice of \emph{antisymmetric tensors} in
$M^{\otimes 2}$, that is,
$$
\wed{M}\simeq \mathsf{A}'_2(M)=\{x\in M^{\otimes 2}\mid x^{\tau}=-x\}
$$
where $\tau:M^{\otimes 2}\to M^{\otimes 2}$ is the switch
$(m\otimes m')^{\tau}=m'\otimes m$; see
\cite[Exerc.~8 on p.~A~III.190]{bourbaki}.
Furthermore,  $\mathsf{A}'_2(M)$ is exactly the kernel of the
canonical map $M^{\otimes 2}\onto\sym{M}$. Hence, we have an exact
sequence of $\cG$-lattices
\begin{equation} \label{E:exact}
0\to\wed{M}\to M^{\otimes 2}\to\sym{M}\to 0 \ .
\end{equation}

\section{Groups of the form $T_{n-1} \sdp \cG$ and lattices}
\label{sect2.5}

\subsection{Notations}\label{S:notations}
In this section we shall focus on the following situation.
Let $T_{n-1} = (\bbG_m)^n/\Delta$ be the (diagonal) maximal
torus of $\PGLn$; here $\Delta \simeq \bbG_m$ diagonally embedded in
$(\bbG_m)^n$.  Recall that $\Sym_n$ acts on $T_{n-1}$ by permuting
the $n$ copies of $\bbG_m$ and that the normalizer $N(T_{n-1})$ of $T_{n-1}$ in $\PGL_n$
is isomorphic to $T_{n-1} \sdp \Sym_n$. We shall be interested in subgroups
of $N(T_{n-1})$ of the form $T_{n-1} \sdp \cG$, 
where $\cG$ is a subgroup of $\Sym_n$. These groups have two properties that
will be important to us in the sequel:
(i) $T_{n-1}\sdp\cG$-varieties and their compressions can 
be constructed from $\cG$-lattices and (ii) certain
$(\PGLn, T_{n-1} \sdp \cG)$-sections will naturally come up in the proofs
of Theorems~\ref{thm1} and~\ref{thm2}. In this section we will focus on
the relationship between $T_{n-1} \sdp \cG$-varieties and $\cG$-lattices.

\subsection{$T_{n-1}\sdp\cG$-varieties and $\cG$-lattices}\label{S:action}
Suppose we are given a morphism
$$
f \colon M \lra A_{n-1}
$$
of $\cG$-lattices, where 
$A_{n-1}$ is the root lattice defined in \eqref{e.A}.
Note that $A_{n-1} \simeq X_*(T_{n-1})$ as an $\Sym_n$-lattice (and hence as $\cG$-lattice), 
where $X_*(T_{n-1})$ is the lattice of characters of $T_{n-1}$. We will always identify
$A_{n-1}$ with $X_*(T_{n-1})$.

We will now associate to $f$ a $T_{n-1}\sdp\cG$-variety $X_f$ as follows. 
Let $\cG$ act on $k[M]$ as usual and define a $T_{n-1}$-action on $k[M]$ by putting
$$
t(m) = f(m)(t)\cdot m \qquad (t\in T_{n-1},\; m\in M)\ .
$$
One easily checks that, by $k$-linear 
extension of this rule, one obtains a well-defined action of $T_{n-1}$ by automorphism on
$k[M]$. Moreover, for $t\in T_{n-1}$, $g\in\cG$ and $m\in M$, one calculates
$$
t(gm)=f(gm)(t)\cdot gm=[gf(m)](t)\cdot gm=f(m)(t^g)\cdot gm=[gt^g](m)\ ;
$$
so the actions of $\cG$ and $T_{n-1}$ combine to yield a locally finite action of
$T_{n-1}\sdp\cG$ on $k[M]$ and thus an algebraic action on $X_f=\Spec k[M]$. 

\begin{lem} \label{lem2.2-1}
The $T_{n-1}\sdp\cG$-variety $X_f$ is a generically free if and only if
\begin{enumerate}
\item
$f$ is surjective and
\item
$\Ker(f)$ is a faithful $\cG$-lattice.
\end{enumerate}
\end{lem}

\begin{proof} Condition (a) is equivalent to saying that the $T_{n-1}$-action on
$X_f$ is generically free; cf., e.g.,~\cite[Theorem 3.2.5]{ov}.
To interpret condition (b) geometrically,
note that $k(X_f)=k(M)$ and $k(X_f/T_{n-1})=k(M)^{T_{n-1}}=k(\Ker(f))$.
Thus condition (b) holds iff $\cG$ acts faithfully on $X_f/T_{n-1}$ or, 
equivalently, iff the $\cG$-action on $X_f/T_{n-1}$ is and
generically free (the two notions coincide
for finite groups). The desired conclusion now follows from Lemma~\ref{lem2.1}.
\end{proof}

\subsection{Compressions and $\cG$-lattices}
In the sequel we will only be interested in generically free 
$T_{n-1}\sdp\cG$-varieties. In particular, we will assume that $f$ is surjective and
expand it into an exact sequence of $\cG$-lattices:
\begin{equation} \label{e.seq}
0 \lra K = \Ker(f) \lra M \stackrel{f}{\lra} A_{n-1} \lra 0 \, ,
\end{equation}
where $K$ is a faithful $\cG$-lattice.
For future reference, we extract the following equality from the proof of
Lemma~\ref{lem2.2-1}:
\begin{equation}\label{E:xfinvar}
k(X_f/ T_{n-1} \sdp \cG) = k(X_f)^{T_{n-1}\sdp\cG} =  [k(M)^{T_{n-1}}]^{\cG}= k(K)^{\cG}\ .
\end{equation}
We can now obtain information about $T_{n-1}\sdp\cG$-compressions of
$X_f$ by studying this sequence more closely.

\begin{lem} \label{lem.diagr}
Suppose that the exact sequence~\eqref{e.seq} 
extends to a commutative diagram
 $$
 \xymatrix{
 0 \ar[r] & K  \ar[r] & M \ar[r]^-{f} & A_{n-1} \ar[r] & 0 \\
 0 \ar[r] & K_0 \ar[r] \ar@{{ _{(}}->}[u] & M_0 \ar[r]^-{f_0} \ar[u] 
 &  A_{n-1} \ar@{=}[u]\ar[r] & 0
 }
$$
of $\cG$-lattices, where $K_0$ is faithful, and the vertical
map $M_0 \lra M$ is injective.
Then there exists a $T_{n-1}\sdp\cG$-compression
$X_f \brokrarr X_{f_0}$.
\end{lem}

\begin{proof} Since $k(M_0)=k(X_{f_0})$ is an $T_{n-1}\sdp\cG$-invariant
subfield of $k(M) = k(X_f)$, it 
defines a dominant $T_{n-1} \sdp \cG$-equivariant map
$X_f \brokrarr X_{f_0}$.
Furthermore, by Lemma~\ref{lem2.2-1}, the $T_{n-1} \sdp \cG$-action
on both $X_f$ and $X_{f_0}$ is generically free. Thus, the rational map 
$X_f \brokrarr X_{f_0}$ we have constructed is a $T_{n-1}\sdp\cG$-compression.
\end{proof}

\subsection{Linearization and essential dimension}\label{S:linearization}
If the $\cG$-lattice $M$ in Section~\ref{S:action} is a permutation lattice
then the $T_{n-1}\sdp\cG$-variety $X_f$ is birationally linearizable.
To see this, fix a $\bbZ$-basis $m_1,\dots,m_r$ of $M$ that is permuted by $\cG$.
Clearly, $k(X_f)=k(M)=k(m_1,\dots,m_r)$. Thus, putting 
$V_f=\sum_i k m_i$ we obtain a $T_{n-1}\sdp\cG$-invariant $k$-subspace of $k(X_f)$
with $k(X_f)=k(V_f)$.

A similar argument goes through if
$M$ to be \emph{permutation projective}, i.e., $M$ is a
direct summand of a permutation $\cG$-lattice.

\begin{lem}\label{L:ess}
If $M$ is permutation projective in \eqref{e.seq} then there is a 
$T_{n-1}\sdp\cG$-compression $V\brokrarr X_f$ with $V$ 
a generically free linear $T_{n-1}\sdp\cG$-variety.
In particular, $\ed(T_{n-1}\sdp\cG) \leq \rank K$.
\end{lem}

\begin{proof}
Suppose $M\oplus N=P$, where $P$ is a permutation $\cG$-lattice.  Then the
sequence~\eqref{e.seq} embeds in the obvious fashion in an exact sequence
$0\to K\oplus N\to P=M\oplus N\to A_{n-1}\to 0$. In view of the foregoing
and Lemma~\ref{lem.diagr}, this proves the first assertion.

To complete the proof, recall that the essential dimension of an algebraic 
group $G$ is defined as the smallest possible value $\dim(X/G)$, where
$X$ is a generically free $G$-variety so that there is a $G$-compression
$V\brokrarr X$ with $V$ a generically free linear $G$-variety; see
\cite[Definition 3.5]{reichstein1}. For $G=T_{n-1}\sdp\cG$ and $X=X_f$, 
we have $\dim(X/G)=\rank K$ by~\eqref{E:xfinvar}.
\end{proof}

\section{Proof of Theorem~\ref{thm1}}
\label{sect3}

\subsection{Reduction to a lattice-theoretic problem}
The universal division algebra $\UD(n)$ is represented by
a class $c \in H^1(Z(n), \PGLn)$.
We can write $\UD(n) = D \otimes_{Z(D)}
Z(n)$ if and only if $c$ lies in the image of the natural map
\begin{equation} \label{e3.1}
H^1(Z(D), \PGLn) \lra H^1(Z(n), \PGLn)
\end{equation}
Recall that for any finitely generated field extension $L/k$,
an element of $\alpha \in H^1(L, \PGLn)$ may also be interpreted
as a {\em $\PGLn$-torsor}, i.e., a generically free
$\PGLn$-variety $X_{\alpha}$
such that $k(X_{\alpha})^{\PGLn} = L$.
Moreover, $X_{\alpha}$ is uniquely determined
(up to birational isomorphism of $\PGLn$-varieties), and
the central simple algebra
corresponding to $\alpha$ can be recovered as the algebra of
$\PGL_n$-equivariant rational maps $X_{\alpha} \brokrarr \Mn$.
In particular, $X_c = \Mn \times \Mn$, where
$\PGLn$ acts on $\Mn \times \Mn$ by simultaneous conjugation;
to say that $c$ lies in the image of the map \eqref{e3.1} is equivalent
to saying that there exists a $\PGLn$-compression
$\Mn \times \Mn \brokrarr X'$ such that $k(X')^{\PGLn} = Z(D)$.
For a more detailed discussion of these facts and further references,
see~\cite[Section 3]{ry}.

Denote the linear subspace of $\Mn$ consisting of diagonal matrices by
$D_n$.  It is easy to see that $D_n \times \Mn$
is a $(\PGLn, N(T_{n-1}))$-section of $\Mn \times \Mn$,
where $N(T_{n-1})=T_{n-1} \sdp \Sym_n$ is the normalizer of the maximal
torus $T_{n-1} = (\bbG_m)^n/\Delta$ in $\PGLn$, as in the previous section.
(See Lemma~\ref{lem5.1} below for a more general fact.)
Thus, in view of Lemma~\ref{lem2.2}, we have the following

\begin{reduction} \label{red3.1}
In order to prove Theorem~\ref{thm1} it is enough to
show that there exists an $N(T_{n-1})$-compression
\begin{equation} \label{e3.3}
 D_n \times \Mn \brokrarr X
\end{equation}
such that
\begin{enumerate}
\item[(i)]
$\dim(X) = \frac{1}{2}(n-1)(n-2) + n-1$ or equivalently,
$\dim(X/N(T_{n-1})) = \trdeg_k k(X)^{N(T_{n-1})} = \frac{1}{2}(n-1)(n-2)$.
\item[(ii)]
$Z(n) = k(\Mn \times \Mn)^{\PGLn}
= k(D_n \times \Mn)^{N(T_{n-1})}$ is purely transcendental over $k(X)^{N(T_{n-1})}$.
\end{enumerate}
\end{reduction}


Our construction of the compression~\eqref{e3.3} will be based
on Lemma~\ref{lem.diagr}. In order to apply this lemma, we need to write
the linear $N(T_{n-1})$-variety $D_n \times \Mn$ birationally in
the form $X_f$, where $f$ is as in~\eqref{e.seq}.
Let $x_i$ and $y_{rs}$ be the standard coordinates on $D_n$ and $\Mn$
respectively. In this coordinate system, the $\Sym_n$-action on 
$D_n \times \Mn$ is given by
$$ 
\sigma(x_i)=x_{\sigma(i)}\quad\text{and}\quad\sigma(y_{rs})=
y_{\sigma(r)\sigma(s)}\; .  
$$
Thus, monomials in these coordinates and their inverses form an $\Sym_n$-lattice
isomorphic to $M = U_n \oplus U_n^{\otimes 2}$, where
$U_n = \bbZ^n$ be the standard permutation $\Sym_n$-lattice.
Moreover, an element 
$t = (t_1, \dots, t_n)$ of $T_{n-1} = (\bbG_m)^n/\Delta$ 
acts on monomials in $x_i$, $y_{rs}$ by characters determined 
(multiplicatively) by
$$
t(x_i) = x_i \; \, \text{and} \; \, t(y_{rs}) = t_r t_s^{-1} y_{rs} \; . 
$$
Denoting the standard basis of $U_n$ by $b_1, \dots, b_n$ 
and defining $f \colon M=U_n \oplus U_n^{\otimes 2} \lra A_{n-1}$
by $f(b_i, b_r \otimes b_s)=b_r - b_s$, the above formulas
give exactly the action of $N(T_{n-1})=T_{n-1}\sdp\Sym_n$ on $X_f$ as described in 
(\ref{S:action}).
The exact sequence~\eqref{e.seq}
for this $f$ is the Formanek -- Procesi exact sequence
\begin{equation} \label{e3.4}
 0 \lra K = \Ker (f) \lra U_n \oplus U_n^{\otimes 2} \stackrel{f}{\lra}
A_{n-1} \lra 0\, ;
\end{equation}
see~\cite{formanek-3x3}.
Note that 
$$
K \simeq U_n\oplus U_n\oplus A_{n-1}^{\otimes 2}\; .
$$
Here, the first copy of $U_n$
is mapped identically onto the first summand of $U_n \oplus U_n^{\otimes 2}$,
the second $U_n$ corresponds to the sublattice of $U_n^{\otimes 2}$ consisting of
the monomials in $y_{ii}\in K$, and
$A_{n-1}^{\otimes 2}$ describes the kernel of $f$, restricted to the sublattice
$\langle y_{rs}\mid r\neq s\rangle\simeq U_n\otimes A_{n-1}\subset U_n^{\otimes 2}$ ;
cf.~\cite[p.~3573]{beneish}.

We now want to apply Lemma~\ref{lem.diagr} to the above sequence, with $K_0 =
\wed{A_{n-1}}$.
Recall that $\wed{A_{n-1}}$ may be viewed as a sublattice $A_{n-1}^{\otimes 2}$;
see~\eqref{E:exact}. Let
\begin{equation} \label{e.phi}
\varphi \colon \wed{A_{n-1}} \hookrightarrow
U_n \oplus U_n \oplus A_{n-1}^{\otimes 2}
\end{equation}
be the natural embedding of $\wed{A_{n-1}}$ into the third component
of $U_n \oplus U_n \oplus A_{n-1}^{\otimes 2}$. 
We remark that for $n > 3$,
$\Sym_n$ acts faithfully on $\wed{A_{n-1}}$; in fact,
$\wed{A_{n-1}} \otimes_{\bbZ} \bbQ$
is the irreducible $S_n$-representation
corresponding to the partition $(n-2,1^2)$
of $n$; see, e.g.,~\cite[Exerc.~4.6]{fulton}.

Combining Reduction~\ref{red3.1} with Lemma~\ref{lem.diagr}, we obtain:

\begin{reduction} \label{red3.2} Theorem~\ref{thm1} follows from
Propositions~\ref{prop3.02} and~\ref{prop3.1} stated below.
\end{reduction}

\begin{prop} \label{prop3.02}
For odd $n \geq 5$,
$k(U_n \oplus U_n\oplus A_{n-1}^{\otimes 2})\simeq k(\wed{A_{n-1}})(y_1,\ldots,y_r)$
as $\Sym_n$-fields,
where the elements $y_i$ are $\Sym_n$-invariant and transcendental
over $k(\wed{A_{n-1}})$.
In particular, $k(U_n \oplus U_n\oplus A_{n-1}^{\otimes 2})^{\Sym_n}$ is rational over
$k(\wed{A_{n-1}})^{\Sym_n}$.
\end{prop}

\begin{prop} \label{prop3.1}
For odd $n$, there exists a commutative diagram of $\Sym_n$-lattices:
$$
\xymatrix{
0 \ar[r] & U_n\oplus U_n\oplus A_{n-1}^{\otimes 2}  \ar[r]
& U_n\oplus U_n^{\otimes 2} \ar[r]^-{f} & A_{n-1} \ar[r] & 0 \\
0 \ar[r] & \wed{A_{n-1}} \ar[r] \ar@{_{(}->}[u]^-{\varphi} & L \ar[r]^-{f_0} \ar[u]
&  A_{n-1} \ar@{=}[u]\ar[r] & 0
}
$$
Here the first row is the Formanek -- Procesi sequence~\eqref{e3.4}.
\end{prop}
Indeed, Proposition~\ref{prop3.1} in conjunction 
with Lemma~\ref{lem.diagr} yields
an $N(T_{n-1})$-compression 
$$
D_n \times \Mn \stackrel{\simeq}{\brokrarr} X_f \brokrarr X=X_{f_0}\;,
$$
and formula~\eqref{E:xfinvar} further implies that 
$k(D_n\times M_n)^{N(T_{n-1})}=k(X_f)^{N(T_{n-1})}=k(U_n\oplus U_n\oplus A_{n-1}^{\otimes 2})^{\Sym_n}$
and $k(X)^{N(T_{n-1})}=k(\wed{A_{n-1}})^{\Sym_n}$. Thus, condition (i) in Reduction~\ref{red3.1}
is clearly satisfied and Proposition~\ref{prop3.02} ensures that (ii) holds as well.

\subsection{Solution of the lattice-theoretic problem}
Our proofs of Propositions~\ref{prop3.02} and~\ref{prop3.1}
will be based on the following result from~\cite[Section 3.5]{ll}.
If $\cG$ is a finite group, $\cH$ is a subgroup of $\cG$ and $M$
a $\bbZ[\cH]$-module then $M\ind{\cH}{\cG}=\bbZ[\cG]\otimes_{\bbZ[\cH]}M$
will denote the induced $\bbZ[\cG]$-module.

\begin{prop}\label{P:ll}
For odd $n$, there is an isomorphism of $\Sym_n$-lattices
$$
\sym{A_{n-1}}\oplus U_n\oplus\bbZ
\simeq \bbZ\ind{\Sym_{n-2}\times\Sym_2}{\Sym_n} \oplus U_n\oplus\bbZ\ ,
$$
where $\bbZ$ has
the trivial $\Sym_n$-action. In particular, $\sym{A_{n-1}}\oplus U_n\oplus\bbZ$
is a permutation lattice.
\end{prop}

\subsubsection*{Proof of Proposition~\ref{prop3.02}}

First, applying Proposition~\ref{P:lenstra} to the obvious sequence
$0\to U_n\oplus A_{n-1}^{\otimes 2}
\to U_n \oplus U_n\oplus A_{n-1}^{\otimes 2}\to U_n\to 0$, we see that
$$
k(U_n \oplus U_n\oplus A_{n-1}^{\otimes 2})
\simeq k(U_n\oplus A_{n-1}^{\otimes 2})(t_1,\ldots,t_n)
$$
as $\Sym_n$-fields.

Next, sequence \eqref{E:exact} for $M=A_{n-1}$ combined with
Proposition~\ref{P:ll} gives rise to an exact sequence of $\Sym_n$-lattices
$$
0\to\wed{A_{n-1}}\to A_{n-1}^{\otimes 2}\oplus U_n\oplus \bbZ \to P\to 0 \ ,
$$
where $P=\sym{A_{n-1}}\oplus U_n\oplus\bbZ$ is permutation.
Applying Proposition~\ref{P:lenstra} to this sequence, we deduce that
$$
k(A_{n-1}^{\otimes 2}\oplus U_n\oplus\bbZ)\simeq k(\wed{A_{n-1}})(x_1,\ldots,x_m)
$$
as $\Sym_n$-fields.

Finally,
\begin{eqnarray*}
k(U_n \oplus U_n\oplus A_{n-1}^{\otimes 2}) &\simeq&
k(U_n\oplus A_{n-1}^{\otimes 2})(t_1,\ldots,t_n)\\
&=& k(A_{n-1}^{\otimes 2}\oplus U_n\oplus\bbZ)(t_1,\ldots,t_{n-1})\\
&\simeq& k(\wed{A_{n-1}})(x_1,\ldots,x_m,t_1,\ldots,t_{n-1})
\end{eqnarray*}
as $\Sym_n$-fields, which proves the first assertion of Proposition~\ref{prop3.02}.
The second assertion is an immediate consequence of the first.
\qed

\subsubsection*{Proof of Proposition~\ref{prop3.1}} Recall that the embedding
$\varphi$ of~\eqref{e.phi} is defined as the composition
$$ \varphi \colon \wed{A_{n-1}} \stackrel{\psi}{\hookrightarrow} A_{n-1}^{\otimes 2}
\hookrightarrow U_n\oplus U_n\oplus A_{n-1}^{\otimes 2}\ , $$
where $\psi$ is the injection from~\eqref{E:exact} (with $M=A_{n-1}$) and the second map
identifies $A_{n-1}^{\otimes 2}$ with the third component of
$U_n \oplus U_n \oplus A_{n-1}^{\otimes 2}$.
We aim to show that $\varphi$ together with sequence \eqref{e3.4} will give
rise to a commutative diagram as in the statement of
Proposition~\ref{prop3.1}. In other words, our goal is
to show that the class in
$\Ext_{\bbZ[\Sym_n]}(A_{n-1},U_n\oplus U_n\oplus A_{n-1}^{\otimes 2})$
corresponding to the extension \eqref{e3.4} belongs to the image of the map
$$
\varphi_\ast \colon \Ext_{\bbZ[\Sym_n]}(A_{n-1},\wed{A_{n-1}}) \to
\Ext_{\bbZ[\Sym_n]}(A_{n-1},U_n\oplus U_n\oplus A_{n-1}^{\otimes 2}) \ .
$$
In fact, we will prove:

\begin{lem}
For odd $n$, the map $\varphi_\ast$ is surjective.
\end{lem}

\begin{proof}
We will tacitly use the following standard facts from homological algebra, valid for
any finite group $\cG$:
\begin{itemize}
\item If $V$ is a $\cG$-module and $M$ an $\cH$-module for some subgroup $\cH\le\cG$
  then $\Ext^\ast_{\bbZ[\cG]}(V,M\ind{\cH}{\cG}) \simeq \Ext^\ast_{\bbZ[\cH]}(V\big|_{\cH},M)$; see
  \cite[Prop. IV.12.3]{hilton} and \cite[Prop. III.5.9]{brown}. For $V=\bbZ$, the trivial
  $\cG$-module, this isomorphism is the ``Shapiro isomorphism"
  $H^\ast(\cG,M\ind{\cH}{\cG})\simeq H^\ast(\cH,M)$. In case, $M$ is actually a $\cG$-module,
  the
  restriction map $\Res_{\cH}^{\cG} \colon H^\ast(\cG,M)\to H^\ast(\cH,M)$ factors through the
  Shapiro isomorphism:
$$
\Res_{\cH}^{\cG}: H^\ast(\cG,M)\stackrel{\mu_\ast}{\to} H^\ast(\cG,M\ind{\cH}{\cG})
\isomoto H^\ast(\cH,M)\ ,
$$
  where $\mu\colon M\to M\ind{\cH}{\cG}$ sends $m\mapsto \sum_{g\in \cG/\cH} g\otimes g^{-1}m$;
  see~\cite[p. 81]{brown}.
\item For any $\cG$-lattices $V$ and $W$,
  $\Ext^\ast_{\bbZ[\cG]}(V,W) \simeq H^\ast(\cG,V^*\otimes W)$,
  where $\otimes=\otimes_{\bbZ}$ and $V^*=\Hom_{\bbZ}(V,\bbZ)$ is the dual $\cG$-lattice; see
  \cite[Prop. III.2.2]{brown}.
\item If $V$ and $W$ are both permutation $\cG$-lattices then
  $\Ext_{\bbZ[\cG]}(V,W)=0$; cf.~\cite[Propositions 1.1, 1.2]{lenstra}.
\end{itemize}
Armed with these facts, we proceed as follows. First,
$\Ext_{\bbZ[\Sym_n]}(A_{n-1},U_n) = 0$, because $U_n \simeq \bbZ\ind{\Sym_{n-1}}{\Sym_n}$
and $A_{n-1}\big|_{\Sym_{n-1}}\simeq U_{n-1}$. Therefore, it suffices to show that the map
$$
\psi_\ast \colon \Ext_{\bbZ[\Sym_n]}(A_{n-1},\wed{A_{n-1}}) \to
\Ext_{\bbZ[\Sym_n]}(A_{n-1},A_{n-1}^{\otimes 2})
$$
is surjective. But the extension \eqref{E:exact} (for $M=A_{n-1}$) gives rise to an
exact sequence
$$
\Ext_{\bbZ[\Sym_n]}(A_{n-1},\wed{A_{n-1}}) \stackrel{\psi_\ast}{\to}
\Ext_{\bbZ[\Sym_n]}(A_{n-1},A_{n-1}^{\otimes 2}) \to
\Ext_{\bbZ[\Sym_n]}(A_{n-1},\sym{A_{n-1}}) \ .
$$
Therefore, it suffices to prove:
$$
\Ext_{\bbZ[\Sym_n]}(A_{n-1},\sym{A_{n-1}}) = 0\quad\text{for odd $n$.}
$$
For this, we use the isomorphism $\sym{A_{n-1}}\oplus U_n\oplus\bbZ
\simeq \bbZ\ind{\cG}{\Sym_n} \oplus U_n\oplus\bbZ$ of Proposition~\ref{P:ll}, where
we have put $\cG=\Sym_{n-2}\times\Sym_2$ for simplicity. This isomorphism
entails
\begin{eqnarray*}
\Ext_{\bbZ[\Sym_n]}(A_{n-1},\sym{A_{n-1}}) &\simeq&
\Ext_{\bbZ[\Sym_n]}(A_{n-1},\bbZ\ind{\cG}{\Sym_n})\\
&\simeq& \Ext_{\bbZ[\cG]}(A_{n-1}\big|_{\cG},\bbZ) \\
&\simeq& H^1(\cG,A_{n-1}^\ast) \ .
\end{eqnarray*}
Dualizing the augmentation sequence $0\to A_{n-1}\to U_n \stackrel{\epsilon}{\to}\bbZ
\to 0$ we obtain an exact sequence
$0\to\bbZ\stackrel{\epsilon^\ast}{\to} U_n\to A_{n-1}^\ast\to 0$, where
$\epsilon^\ast(1)=\sum_ie_i$\;, the sum of the natural basis elements of $U_n$.
This sequence, viewed as exact sequence of $\cG$-lattices,
in turn yields an exact sequence
$$
H^1(\cG,U_n)=0 \to H^1(\cG,A_{n-1}^\ast) \to H^2(\cG,\bbZ) \to H^2(\cG,U_n) \ .
$$
Thus, it suffices to show that $H^2(\cG,\bbZ) \to H^2(\cG,U_n)$ is injective.
As a $\cG$-module, $U_n=V\oplus W$, where
$V=\bigoplus_{i=1}^{n-2}\bbZ e_i \simeq \bbZ\ind{\Sym_{n-3}\times\Sym_2}{\cG}$
and $W=\bbZ e_{n-1}\oplus \bbZ e_n \simeq \bbZ\ind{\Sym_{n-2}}{\cG}$. Therefore, the
Shapiro isomorphism gives
$$
H^2(\cG,U_n) \simeq H^2(\Sym_{n-3}\times\Sym_2,\bbZ) \oplus H^2(\Sym_{n-2},\bbZ)
$$
and the map $H^2(\cG,\bbZ) \to H^2(\cG,U_n)$ becomes the restriction map
$$
\Res_{\Sym_{n-3}\times\Sym_2}^{\cG} \times \Res_{\Sym_{n-2}}^{\cG}\colon
H^2(\cG,\bbZ) \to H^2(\Sym_{n-3}\times\Sym_2,\bbZ) \oplus H^2(\Sym_{n-2},\bbZ) \ .
$$
This map is indeed injective, as is easily seen by
identifying $H^2(\cG,\bbZ)$ in the usual fashion with $\Hom(\cG,\bbQ/\bbZ)$,
and similarly for the subgroups $\Sym_{n-3}\times\Sym_2$ and $\Sym_{n-2}$.
This finishes the proof of the Lemma, and hence, of Proposition~\ref{prop3.1}
and of Theorem~\ref{thm1}.
\end{proof}

\section{Proof of Theorem~\ref{thm2}}
\label{sect5}

\subsection{General observations}\label{S:observations}
The following notations 
will be used throughout this section: \bigskip
\begin{tabbing}
\hspace*{.35in}\=\hspace{.4in}\=\kill
\> $F$ \> will be a field containing $k$;\\
\> $A$ \> will be a finite-dimensional
	central simple algebra with center $F$;\\
\> $L$ \> will be a strictly maximal commutative subfield of $A$;\\
\> $n$ \> denotes the degree of $A$, so $\dim_FA=n^2$ and $[L:F]=n$;\\
\> $\cG$ \> is the Galois group of the normal closure $L^{norm}$ 
	of $L$ over $F$; \\
\> $T_{n-1}$ \> denotes the diagonal maximal torus of $\PGLn$.\\
\end{tabbing}

\begin{reduction} \label{red5.1} In the course of proving Theorem~\ref{thm2} we may
assume without loss of generality that $F$ is a finitely generated field extension
of $k$.
\end{reduction}

\begin{proof}
Indeed, choose a primitive element $x$ for the extension $L/F$ and complete
$1, x, \dots, x^{n-1}$ to an $F$-basis $e_1, \dots, e_{n^2}$ of $A$, where
$e_i = x^{i-1}$ for $i = 1, \dots, n$. Let $c_{rs}^t$ be the structure
constants for $A$ in this basis, i.e.,
\begin{equation} \label{e5.1}
 e_r e_s = \sum_{t=1}^{n^2} c_{rs}^t e_t
\end{equation}
for every $r, s = 1, \dots, n^2$. Then $A = A_0 \otimes_{F_0} F$,
where $F_0 = k(c_{rs}^t)$ and $A_0$ is the $n^2$-dimensional $F_0$-algebra
spanned by $e_1, \dots, e_{n^2}$ with multiplication given
by~\eqref{e5.1}. Moreover, $A_0$ is a central simple algebra of degree $n$ with center $F_0$
and $A_0$ contains the strictly maximal
subfield $L_0 = F_0(x)$ whose normal closure $L_0^{norm}$ has Galois group $\cG$
over $F_0$. Clearly, $\tau(A) \leq \tau(A_0)$. Thus, after replacing $A$ by $A_0$
we may assume that $F$ is finitely generated extension of $k$.
\end{proof}

Next we pass from central simple algebras
to generically free $\PGL_n$-varieties, as we did at the beginning
of Section~\ref{sect3}. Recall that a central simple algebra $A$ of degree
$n$ with center $F$ defines a class $H^1(F, \PGLn)$.
This class, in turn, gives rise to a $\PGLn$-torsor $X_{A}$ over $F$.
Since we are assuming $F$ is a finitely generated extension of $k$, $X_A$ is a
generically free $\PGLn$-variety; moreover, $F = k(X_A/\PGLn)$ and
$A$ can be recovered from $X_A$
as the algebra of $\PGL_n$-equivariant rational maps $X_A \brokrarr \Mn$;
cf.~\cite[Section 3]{ry}.

The algebras $A$ we are concerned with in the context Theorem~\ref{thm2} are
of a special form: they have a maximal subfield $L/F$ such
that $\Galois(L^{norm}/F)$ is the given group $\cG$. 
We would like to know how this extra structure
is reflected in the geometry of the $\PGLn$-variety $X_A$. The following lemma
gives a partial answer. We continue to let $T_{n-1}$ denote the diagonal maximal torus
of $\PGLn$, as in Section~\ref{S:notations}.

\begin{lem} \label{lem5.1} 
$X_A$ has a $(\PGLn, T_{n-1}\sdp\cG)$-section.
\end{lem}

Note that the group $\cG$ comes with a natural permutation representation
of $\cG$ on the $n$ embeddings $L \hookrightarrow L^{norm}$ over $F$. 
This permutation representation gives 
an embedding $\alpha \colon \cG \hookrightarrow \Sym_n$ that we use
to define the semidirect product $T_{n-1} \sdp \cG$. Note that $\alpha$ is only
defined up to an inner automorphism of $\Sym_n$, since the $n$ embedding
$L \hookrightarrow L^{norm}$ are not naturally in a 1-1 correspondence
with $\{ 1, \dots, n \}$.  Relabeling these embeddings (or, equivalently,
reordering the roots of a defining polynomial for $L/F$) will
cause $T_{n-1} \sdp \cG$ to be replaced by a conjugate subgroup of $\PGLn$;
the corresponding section will be translates of each other by elements
of $\Sym_n \subset \PGLn$. Thus it is sufficient to verify that such
a section exists for a particular numbering of the embeddings
$L \hookrightarrow L^{norm}$.

\begin{proof}
Suppose $L = F(r)$ for some $r \in L$ and $r = r_1, r_2, \dots, r_n$
are the conjugates of $r$ in $L^{norm}$. As we remarked above, the
order of the roots is not intrinsic; however, we choose it at this point,
and it will not be changed in the sequel.  The group $\cG$ permutes 
$r_1, \dots, r_n$ transitively via a permutation representation 
$\alpha \colon \cG \hookrightarrow \Sym_n$.

Now let $x_1, \dots, x_n$ be commuting independent variables over $k$.
The symmetric group acts on the polynomial ring $k[x_1, \dots, x_n]$ by
permuting these variables; composing this action with $\alpha$, we obtain
a (permutation) action of $\cG$ on $k[x_1, \dots, x_n]$. Note that for every
$f(x_1, \dots, x_n) \in k[x_1, \dots, x_n]^{\cG}$, we have $f(r_1, \dots, r_n)
\in F$; we shall write this element of $F$ as $f_r$.

Recall that $A$ is the algebra of $\PGLn$-equivariant rational maps
$X_A \brokrarr \Mn$.  We claim that
\[ S = \left\{ x \in X_A \quad : \quad  \begin{array}{c}
r(x) = \diag(\lambda_1, \dots, \lambda_n)
\quad \text{\rm is a diagonal matrix} \\
\text{\rm and} \\
f(\lambda_1, \dots, \lambda_n) = f_r(x) \quad \text{\rm for every
$f \in k[x_1, \dots, x_n]^{\cG}$} \end{array} \right\} \]
is an $(\PGLn, T_{n-1}\sdp\cG)$-section $S$ of $X_A$. Note that
$S$ is a $T_{n-1}\sdp\cG$-invariant subvariety of $X_A$: indeed,
$T_{n-1}$ acts trivially on the set of diagonal matrices, and $f_r \in F$ 
is a $\PGL_n$-invariant rational function on $X_A$. Thus we need to show
that  

\medskip
(i) $\PGL_n x$ intersects $S$
for $x$ in general position in $X_A$ and

\smallskip
(ii) $gs \in S  \quad \Longrightarrow  \quad g \in T_{n-1}\sdp\cG$
for $s$ in general position in $S$;

\smallskip
\noindent
cf. Section~\ref{sect2.2}.

\smallskip
Let $\pi \colon X_A \brokrarr X_A/\PGLn$ be the rational quotient map.
Recall that $k(X_A/\PGLn)= k(X_A)^{\PGLn} = F$. 
Suppose $p(t) = t^n + a_1t^{n-1} + \dots + a_n$
is the minimal polynomial of $r \in L$ over $F$. Note $a_1, \dots, a_n \in F$
are $\PGLn$-invariant rational functions on $X_A$; in particular, 
for $x \in X_A$ in general position, the matrix $r(x) \in \Mn$
satisfies the polynomial $p_x(t) =
t^n + a_1(x)t^{n-1} + \dots + a_n(x) \in k[t]$. 
Since $p(t)$ is an irreducible polynomial over $F$, 
its discriminant $\delta$ is a non-zero element of $F$, i.e., a non-zero
$\PGLn$-invariant rational function $X_A$. This means that for $x \in X_A$ in general
position (i.e., away from the zero locus of $\delta$ and the indeterminacy
locus of $r$), the $n \times n$-matrix $r(x)$ has distinct eigenvalues. 
We conclude that in this case  $p_x(t)$ is the characteristic polynomial 
for $r(x)$; in particular, the eigenvalues of $r(x)$ are precisely the roots of
$p_x(t)$.

To see what these eigenvalues are more explicitly,
let $Y \brokrarr X_A/\PGL_n$ be a rational map of varieties induced
by the field extension $L^{norm}/F$. Then $r_1, \dots, r_n \in L^{norm}$
are rational functions on $Y$. Thus we have the following diagram of rational 
maps 
$$
\xymatrix@=1pt@!{
X_A \ar@{-->}[dd]_{\pi} \ar@{-->}[r]^r & \Mn \\
& Y \ar@{-->}[r]^{r_i} \ar@{-->}[dl] & k \\
X_A/\PGLn
}
$$
Suppose $x$ be a point of $X_A$ in general position and $y$ is a point of
$Y$ lying above $\pi(x)$.
Since  $r_1, \dots, r_n$ are distinct elements of $L^{norm} = k(Y)$,
$\lambda_1 = r_1(y), \dots, \lambda_n =
r_n(y)$ are the $n$ distinct roots of $p_x(t)$, i.e., the $n$ distinct 
eigenvalues of the $n \times n$-matrix $r(x)$.

\smallskip
Proof of (i): In view of the above discussion, we may assume without loss
of generality that $r(x)$ is diagonalizable.
In other words, the $\PGL_n$-orbit of $r(x)$ in
$\Mn$ contains the diagonal matrix $\diag(\lambda_1, \dots, \lambda_n)$
or equivalently,
$r(x') =   \diag(\lambda_1, \dots, \lambda_n)$ for some $x' \in
\PGLn x$. It remains to show that $x' \in S$.
Indeed,  for any $f \in k[x_1, \dots, x_n]^{\cG}$, we have
\[ f_r(x') = f(r_1(y), \dots, r_n(y)) = f(\lambda_1, \dots, \lambda_n) \, , \]
as desired. This completes the proof of (i).

Proof of (ii): Let $x$ be a point of $S$ in general position. We may
assume without loss of generality that the eigenvalues $\lambda_i$
of the diagonal matrix $r(x) = \diag(\lambda_1, \dots, \lambda_n)$
are distinct. Suppose $gx \in S$ for some $g \in \PGLn$. Then
$gr(x)g^{-1}$ is again diagonal; hence, $g \in T_{n-1} \sdp\Sym_n$.
In other words, $g = t \sigma$, where $t$ is a diagonal matrix and
$\sigma$ is a permutation matrix; our goal is to show that
$\sigma \in \cG$.  Indeed, by the definition of $S$,
$\sigma$ has the property that $f(\lambda_{\sigma(1)}, \dots,
\lambda_{\sigma(n)}) = f(\lambda_1, \dots, \lambda_n)$ for every
$f \in k[x_1, \dots, x_n]^{\cG}$. 
Since $\cG$-invariant regular functions separate the orbits of the permutation
$\cG$-action on affine $n$-space (this is true for any finite group action on an
affine variety; see, e.g., \cite[Section 0.4]{pv}), the points
$(\lambda_{\sigma(1)}, \dots,
\lambda_{\sigma(n)})$ and $(\lambda_1, \dots, \lambda_n)$ are in
the same $\cG$-orbit for this action. On the other hand, since
$\lambda_1, \dots, \lambda_n$ are distinct, this is only possible
if $\sigma \in\cG$. This completes
the proof of (ii) and thus of Lemma~\ref{lem5.1}.
\end{proof}

\begin{remark} \label{rem.converse}
One can show that the converse of Lemma~\ref{lem5.1} is also true:
if $X_A$ has a $(\PGLn, T_{n-1} \sdp \cG)$-section then $A$ contains 
a strictly maximal subfield $L$ such that $\Galois(L^{norm}/F) = \cG$.
Since this result is not needed in the sequel, we omit the proof.
\end{remark} 

\subsection{Conclusion of the proof}

We are now ready to finish the proof of Theorem~\ref{thm2}. 
Let $A$ be a central simple algebra of degree $n$ and let $X_A$ 
denote the $\PGLn$-variety associated to $A$.
Recall that $\tau(A) = \ed(X_A, \PGLn)$;
see~\cite[Theorem 8.8 and Lemma 9.1]{reichstein1}. Moreover,
if $X$ has a $(\PGLn, H)$-section $S$ then $\ed(X, \PGLn) \leq \ed(S, H)
\leq \ed(H)$; see \cite[Lemma 4.1 and Definition 3.5]{reichstein1}.
Applying these inequalities to the situation described
by Lemma~\ref{lem5.1}, with $H = T_{n-1}\sdp\cG$,
we see that
\begin{equation} \label{e.ineq}
\tau(A) \leq \ed(T_{n-1} \sdp \cG) \, .
\end{equation}
Thus Theorem~\ref{thm2} is a consequence of the following:

\begin{lem} \label{lem5.2}
Suppose $\cG$ is a transitive subgroup of $\Sym_n$ generated by
the subgroup $\cH=\cG\cap\Sym_{n-1}$ together with elements $g_1, \dots, g_r$.
Assume that either $r \ge 2$ or $\cH \neq \{ 1\}$.
Then
$\ed(T_{n-1}\sdp\cG) \leq r|\cG|-n+1$.
\end{lem}

\begin{proof}
By Lemma~\ref{L:ess}, it suffices to construct an exact sequence~\eqref{e.seq}
with $M$ permutation projective and $K$ faithful having $\rank K=r|\cG|-n+1$.
To this end, note that $U_n\simeq\bbZ[\cG/\cH]$ as $\cG$-lattices. Let
$\overline{\phantom{iii}}\colon \bbZ[\cG]\onto \bbZ[\cG/\cH]=U_n$ denote 
the canonical epimorphism; the kernel of $\overline{\phantom{iii}}$ is
$\bbZ[\cG]\omega\cH$, where $\omega\cH$ denotes the augmentation
ideal of $\bbZ[\cG]$; cf.~\cite{passman}.
Then $\sum_i\bbZ[\cG]\overline{(g_i-1)}=A_{n-1}$;
see~\cite[Lemma 3.1.1]{passman}. Therefore, we obtain an epimorphism of 
$\cG$-lattices
$$
f\colon M=\bbZ[\cG]^r\onto A_{n-1}\;,\qquad (\alpha_1,\dots,\alpha_r)\mapsto
\sum_i^r\alpha_i\overline{(g_i-1)}\ .
$$
Put $K=\Ker f$; so $K$ certainly has the required rank. For faithfulness, we may
consider $K\otimes\bbQ$ instead of $K$ and work over the semisimple algebra $\bbQ[\cG]$.
Since $f\otimes\bbQ$ and $\overline{\phantom{iii}}\otimes\bbQ$ are split, we
have $\bbQ[\cG]$-isomorphisms 
$(A_{n-1}\otimes\bbQ) \oplus (K\otimes\bbQ) \simeq \bbQ[\cG]^r$
and 
$(A_{n-1}\otimes\bbQ) \oplus \bbQ \oplus \bbQ[\cG]\omega\cH \simeq \bbQ[\cG]$.
Therefore,
$$
K\otimes\bbQ \simeq \bbQ[\cG]^{r-1} \oplus \bbQ \oplus \bbQ[\cG]\omega\cH\ .
$$
If $r\ge 2$ then $\bbQ[\cG]^{r-1}$ is $\cG$-faithful, and 
if $\cH \neq \{ 1\}$ then $\omega\cH\otimes\bbQ$ is $\cH$-faithful and so 
$\bbQ[\cG]\omega\cH \simeq (\omega\cH\otimes\bbQ)\ind{\cH}{\cG}$ is $\cG$-faithful.
In either case, $K\otimes\bbQ$ is faithful, and hence so is $K$, as desired.
\end{proof}


\section{Algebras of degree four}
\label{sect6}
%
%

Recall that Corollary~\ref{cor3} asserts that $d(4)$ equals $4$ or $5$.
Whether the true value of $d(4)$ is four or five is an open question.
The purpose of this section is to show that this question cannot be resolved 
by the methods of this paper.
 
For the rest of this section we will identify the Klein 4-group 
$\cV = (\bbZ/2\bbZ) \times (\bbZ/2\bbZ)$ with the subgroup 
of $\Sym_4$ generated by $(12)(34)$ and $(13)(24)$. 
Let $A_3$ be the augmentation (or root) lattice of $\Sym_4$, 
restricted to $\cV$; see ~\eqref{e.A}. In other words, 
\begin{equation} \label{e.omega}
A_3\simeq \omega\cV\, , 
\end{equation}
where $\omega \cV$ is the augmentation ideal of the group ring $\bbZ[\cV]$.

We now briefly recall
how we arrived at the bound $d(4) \leq 5$. First of all, since
$\UD(4)$ is a $\cV$-crossed product,
$d(4)  = \tau(\UD(4)) \leq \ed(T_3 \sdp \cV)$; see~\eqref{e.ineq}. Secondly,
Lemma~\ref{lem.diagr} tells us that
$\ed(T_3 \sdp \cV) \leq \rank(K_0)$, for any
commutative diagram 
\begin{equation}\label{diagr}
\xymatrix{ 0 \ar[r] & K  \ar[r] & M \ar[r]^-{f} & A_3 \ar[r] & 0 \\
 0 \ar[r] & K_0 \ar[r] \ar[u]^{\varphi} & M_0 \ar[r]^-{f_0} \ar[u] 
 &  A_3 \ar@{=}[u]\ar[r] & 0 }
\end{equation}
of $\cV$-lattices with $M$ permutation projective, $K_0$ faithful, 
and $\varphi$ is injective.  Finally, in the course of the proof 
of Lemma~\ref{lem5.2} (with $\cG = \cV$ and $r = 2$) we constructed 
a particular diagram~\eqref{diagr} with $M = M_0 = \bbZ[\cV]^2$
and  $\rank(K_0) = 5$. This gave us the bound $d(4) \leq 5$.
The question we will now address is whether or not one can sharpen this
bound by choosing a different diagram~\eqref{diagr}. The following 
proposition shows that the answer is ``no".

\begin{prop} \label{prop6.2} Let~\eqref{diagr} be a commutative diagram
of $\cV = (\bbZ/2\bbZ) \times (\bbZ/2\bbZ)$-lattices with $M$ permutation projective.
Then $K_0$ is faithful and $\rank K_0\ge 5$.
\end{prop}

Note that in the setting of Lemma~\ref{lem.diagr}
we assumed 
that (i) $K_0$ is faithful and (ii) $\varphi$ is injective. Here
we see that (i) is automatic and (ii) is irrelevant for the rank estimate.

\begin{proof}
Since $M$ is permutation projective, we have $H^1(\cH,M)=0=H^{-1}(\cH,M)$ 
for every subgroup $\cH$ of $\cV$. 
(This condition is actually equivalent to $M$ being
permutation projective; see \cite[Proposition 4]{CTS}.) Consequently, 
\eqref{diagr} yields a commutative diagram
$$
\xymatrix{
 0 = H^1(\cH,M)  \ar[r] & H^1(\cH,A_3) \ar[r]^-{\delta} & H^2(\cH,K) \\
 & H^1(\cH,A_3) \ar[r]^-{\delta_0} \ar@{=}[u] & H^2(\cH,K_0) \ar[u]^{\varphi_*}
 }
$$
Thus, $\delta_0$ is mono. Since $H^1(\cV,A_3)\simeq \bbZ/4\bbZ$ and 
$H^1(\cH,A_3)\simeq \bbZ/2\bbZ$ for any nonidentity cyclic subgroup 
$\cH$ of $\cV$ (see \cite[Lemma 4.3]{ll}), we obtain
 \begin{equation}\label{twocohom}
 \bbZ/4\bbZ\hookrightarrow H^2(\cV,K_0)\qquad\text{and}\qquad\text{$H^2(\cH,K_0)\neq 0$
 for all $1\neq\cH\le\cV$.}
 \end{equation}
Similarly, $0 = H^{-1}(\cH,M)$ implies
$H^{-1}(\cH,A_3) \hookrightarrow \widehat{H}^0(\cH,K_0)$. Using the 
identification~\eqref{e.omega}, we have 
$H^{-1}(\cV,A_3)\simeq \omega\cV/(\omega\cV)^2\simeq \bbZ/2\bbZ\oplus\bbZ/2\bbZ$.
Thus:
\begin{equation}\label{zerocohom}
\bbZ/2\bbZ\oplus\bbZ/2\bbZ \hookrightarrow \widehat{H}^0(\cH,K_0)\ . 
\end{equation} 
We will show that \eqref{twocohom} forces $K_0$ to be faithful,
and \eqref{twocohom}  
and \eqref{zerocohom} together imply that $\rank K_0\ge 5$.
The discussion below could be shortened somewhat by a reference 
to~\cite{naz}; however, for the sake of completeness, 
we will give a self-contained argument.

\begin{lem} \label{lem6.4}
Let $L$ be a $\cV$-lattice, $1\neq x\in\cV$ and
$L_{\pm}=\{l\in L\mid xl=\pm l\}$. 
If $L\big|_{\langle x\rangle}=L_{+}\oplus L_{-}$
then $2\cdot H^2(\cV,L)=0$.
\end{lem}

\begin{proof} Since $L_{+}$ and $L_{-}$ are $\cV$-sublattices of $L$,
we may assume $L=L_{+}$ or $L=L_{-}$. Write $\cV=\langle x,y\rangle$. Then
$\langle y\rangle$-sublattices of $L$ are stable under $\cV$. 
Therefore, we may assume that $L$ is indecomposable 
as a $\langle y\rangle$-lattice. This leaves the following possibilities 
for $L$: $\bbZ_{\pm}\ind{\langle x\rangle}{\cV}$ or
$\bbZ_{\lambda}$ for some $\lambda\in\Hom(\cV,\bbZ)$. 
In each case, $2\cdot H^2(\cV,L)=0$ is easy to verify.
\end{proof}

Lemma~\ref{lem6.4}, in combination with the first condition 
in~\eqref{twocohom}, implies that $K_0$ is faithful.
Consequently, $\overline{K_0}=K_0/K_0^{\cV}$ is faithful as well, and so
$\rank\overline{K_0}\ge 2$. In addition, we know by \eqref{zerocohom}
that $\widehat{H}^0(\cV,K_0)=K_0^{\cV}/(\sum_{\cV}v)K_0$ is not cyclic. Hence,
neither is $K_0^{\cV}$, which forces $\rank K_0^{\cV}\ge 2$ and thus
$\rank K_0\ge 4$. Suppose, by way of contradiction, that equality holds 
here, i.e., $\overline{K_0}=K_0/K_0^{\cV}$ and $K_0^{\cV}$ both 
have rank 2. By the well-known classification of finite subgroups 
of $\GL_2(\bbZ)$, the action of $\cV$ on $\overline{K_0}$
is given by either
\begin{description}
\item[diag] the matrices $\begin{pmatrix}1 & \\ & -1\end{pmatrix}$ and 
$\begin{pmatrix}-1 & \\ & 1\end{pmatrix}$;
so $\overline{K_0}\simeq \bbZ_{+,-}\oplus\bbZ_{-,+}$, or
\item[non-diag] the matrices $\begin{pmatrix} & 1\\ 1 & \end{pmatrix}$ and 
$\begin{pmatrix} & -1 \\ -1 & \end{pmatrix}$; 
so $\overline{K_0}\simeq \bbZ_{-}\ind{\cH}{\cV}$ for some cyclic $\cH\le \cV$.
\end{description}
In case \textbf{non-diag}, $H^2(\cV,\overline{K_0})\simeq H^2(\cH,\bbZ_{-})\simeq
\widehat{H}^0(\cH,\bbZ_{-})=0$, and hence $H^2(\cV,K_0^{\cV})$ maps onto $H^2(\cV,K_0)$.
But $H^2(\cV,K_0^{\cV})\simeq H^2(\cV,\bbZ)^2\simeq \Hom(\cV,\bbQ/\bbZ)^2\simeq (\bbZ/2\bbZ)^4$.
Thus $H^2(\cV,K_0^{\cV})$ is annihilated by $2$, and hence 
so is $H^2(\cV,K_0)$, 
contradicting \eqref{twocohom}. Therefore, \textbf{diag} must hold:
$$
\overline{K_0}=K_0/K_0^{\cV} \simeq \bbZ_{+,-}\oplus\bbZ_{-,+}\ .
$$
The action of $\cV$ on $K_0$ is given by matrices
$$
c=\left(
\begin{array}{c|c}
\boldsymbol{1}_{2\times 2} & 
\begin{matrix} \boldsymbol{0} & \boldsymbol{\gamma}  \end{matrix} \\
\hrulefill & \hrulefill \\
& \begin{matrix} 1 & \\ & -1 \end{matrix}
\end{array}
\right)
\qquad\text{and}\qquad
d=\left(
\begin{array}{c|c}
\boldsymbol{1}_{2\times 2} & 
\begin{matrix} \boldsymbol{\delta} & \boldsymbol{0} \end{matrix} \\
\hrulefill & \hrulefill \\
& \begin{matrix} -1 & \\ & 1 \end{matrix}
\end{array}
\right)
$$
with $\boldsymbol{\gamma}, \boldsymbol{\delta}\in M_{2\times 1}(\bbZ)$ and 
$\boldsymbol{0}=\left(\begin{smallmatrix}0 \\ 0\end{smallmatrix}\right)$.
By Lemma~\ref{lem6.4}, $\boldsymbol{\gamma}\neq \boldsymbol{0}$ and
$\boldsymbol{\delta}\neq \boldsymbol{0}$.
Conjugating by a suitable matrix of the form $\left(
\begin{smallmatrix}
\boldsymbol{1}_{2\times 2} & 
\begin{smallmatrix} \boldsymbol{0} & \boldsymbol{\rho} \end{smallmatrix} \\
& \boldsymbol{1}_{2\times 2} 
\end{smallmatrix}
\right)$
we can ensure that the entries of $\boldsymbol{\gamma}$
are $0$ or $1$, 
and similarly for $\boldsymbol{\delta}$.
If $\boldsymbol{\gamma}=\left(\begin{smallmatrix}1 \\ 1\end{smallmatrix}\right)$
or $\left(\begin{smallmatrix}1 \\ 0\end{smallmatrix}\right)$
then conjugating, respectively, by $\left(
\begin{smallmatrix}
\begin{smallmatrix} 1 & 1 \\ & 1\end{smallmatrix} & \\
& \boldsymbol{1}_{2\times 2} 
\end{smallmatrix}
\right)$
or 
$\left(
\begin{smallmatrix}
\begin{smallmatrix}  & 1 \\1 & \end{smallmatrix} & \\
& \boldsymbol{1}_{2\times 2} 
\end{smallmatrix}
\right)$, 
we can replace $\boldsymbol{\gamma}$ by 
$\boldsymbol{\gamma}=\left(\begin{smallmatrix}0 \\ 1\end{smallmatrix}\right)$. 
Thus we may assume that 
$\boldsymbol{\gamma}=\left(\begin{smallmatrix}0 \\ 1\end{smallmatrix}\right)$. 
If $\boldsymbol{\delta}=\left(\begin{smallmatrix}1 \\ 1\end{smallmatrix}\right)$, then
conjugating by $\left(
\begin{smallmatrix}
\begin{smallmatrix} 1 & \\ 1 & 1\end{smallmatrix} & \\
& \boldsymbol{1}_{2\times 2} 
\end{smallmatrix}
\right)$, we replace $\boldsymbol{\delta}$ by 
$\left(\begin{smallmatrix}1 \\ 0\end{smallmatrix}\right)$ without changing $c$.
This leaves us with two cases to consider:
\bigskip

\underline{$\boldsymbol{\delta}=\left(\begin{smallmatrix}0 \\ 1\end{smallmatrix}\right)$}:
Then $c=\left(\begin{smallmatrix}1 & \\ & c'\end{smallmatrix}\right)$ and 
$d=\left(\begin{smallmatrix}1 & \\ & d'\end{smallmatrix}\right)$ with
$c'=\left(\begin{smallmatrix}1 & 0 & 1 \\  0 & 1 & 0 \\ 0 & 0 & -1\end{smallmatrix}\right)$
and 
$d'=\left(\begin{smallmatrix}1 & 1 & 0 \\  0 & -1 & 0 \\ 0 & 0 & 1\end{smallmatrix}\right)$.
Therefore, $K_0\simeq \bbZ\oplus (A_3\otimes_{\bbZ}\bbZ_{\lambda})$, where
$\lambda\colon \cV\to \bbZ$ is the map sending the elements of $\cV$ acting via
$c$ and $d$ both to $-1$. Tensoring the augmentation sequence 
$0\to A_3=\omega\cV\to \bbZ[\cV]\to\bbZ\to 0$ with $\bbZ_{\lambda}$ we obtain an exact sequence
$0\to A_3\otimes\bbZ_{\lambda} \to \bbZ[\cV]\otimes\bbZ_{\lambda}=\bbZ[\cV]
\to\bbZ_{\lambda}\to 0$. This sequence in turn implies that 
$H^2(\cV,A_3\otimes\bbZ_{\lambda})\simeq H^1(\cV,\bbZ_{\lambda})$, 
and the inflation-restriction sequence easily gives 
$H^1(\cV,\bbZ_{\lambda})=\bbZ/2\bbZ$. Thus, 
$H^2(\cV,K_0)=H^2(\cV,\bbZ)\oplus H^2(\cV,A_3\otimes\bbZ_{\lambda})
\simeq (\bbZ/2\bbZ)^3$, contradicting \eqref{twocohom}.
\bigskip

\underline{$\boldsymbol{\delta}=\left(\begin{smallmatrix}1 \\ 0\end{smallmatrix}\right)$}:
In this case, $cd=\left(
\begin{smallmatrix}
\boldsymbol{1}_{2\times 2} & \boldsymbol{1}_{2\times 2} \\
& -\boldsymbol{1}_{2\times 2} 
\end{smallmatrix}
\right)$.
Letting $\cH$ denote the cyclic subgroup of $\cV$ acting via $cd$, we have
$K_0\big|_{\cH}\simeq \bbZ[\cH]^2$. Thus, $H^2(\cH,K_0)=0$, contradicting 
\eqref{twocohom}.
\bigskip

This completes the proof of the proposition.
\end{proof}


%
\begin{ack} The work on this article was started while the authors were
attending the Noncommutative Algebra program at MSRI in the Fall of 1999. 
The authors would like to thank the organizers of this program and
MSRI staff for their hospitality and support.
\end{ack}
%


\end{document}